\documentclass{article}
\usepackage{amssymb}
\begin{document}
\input psbox.tex
\psfordvips

\centerline{\LARGE\bf On the classification of convex lattice}

\medskip
\centerline{\LARGE\bf polytopes (II)}

\bigskip\bigskip \centerline{\bf Chuanming Zong\footnote{This work is supported by 973
Program 2011CB302400, the National Natural Science Foundation of China (Grant No. 11071003),
the Chang Jiang Scholars Program of China and LMAM at Peking University.}}

\bigskip
\noindent {\bf Abstract.} In 1980, V.I. Arnold studied the classification
problem for convex lattice polygons of given area. Since then this problem and its
analogues have been studied by several authors, upper bounds for the numbers of
non-equivalent $d$-dimensional convex lattice polytopes of given volume or fixed number of lattice points
have been achieved. In this paper, by introducing and studying the unimodular groups
acting on convex lattice polytopes, we obtain a lower bound for the number of non-equivalent
$d$-dimensional centrally symmetric convex lattice polytopes of given number of lattice points,
which is essentially tight.

\bigskip
\noindent 2010 Mathematics Subject Classification. 52B20, 52C07

\section*{1. Introduction}

Let $\{ {\bf e}_1, {\bf e}_2, \ldots , {\bf e}_d\}$ be an orthonormal basis of the
$d$-dimensional Euclidean space $\mathbb{E}^d$. A convex lattice polytope in $\mathbb{E}^d$
is the convex hull of a finite subset of the integral lattice $\mathbb{Z}^d$. As usual,
let $P$ denote a $d$-dimensional convex lattice polytope, let $v(P)$ denote the volume
of $P$, let ${\rm int}(P)$ denote the interior of $P$, and let $|P|$ denote the cardinality
of $P\cap \mathbb{Z}^d$. For general references on polytopes and lattice polytopes, we refer
to \cite{bara08}, \cite{barv04}, \cite{grit93}, \cite{grub07} and \cite{zieg95}.

Let $P_1$ and $P_2$ be $d$-dimensional convex lattice polytopes. If there is a unimodular
transformation $\sigma $ ($\mathbb{Z}^d$-preserving affinely linear transformation) satisfying
$P_2=\sigma (P_1),$ then we say $P_1$ and $P_2$ are equivalent. For convenience, we write $P_1 \sim P_2$ for short.
It is easy to see that, if $P_1\sim P_2$ and $P_2\sim P_3$, then we have $P_1\sim P_3$. In
addition, if $P_1\sim P_2$, then we have $v(P_1)=v(P_2)$ and $|P_1|=|P_2|.$

Clearly, the equivalence relation $\sim$ divides convex lattice polytopes into different classes. Using triangulations, it can be easily shown that $d!\cdot v(P)\in \mathbb{Z}$
holds for any $d$-dimensional convex lattice polytope $P$. Let $v(d,m)$ denote the number of
different classes of the $d$-dimensional convex lattice polytopes $P$
with $v(P)=m/d!$, where both $d$ and $m$ are positive integers.

Let $f(d,m)$ and $g(d,m)$ be functions of positive integers $d$ and $m$. In this paper $f(d,m)\ll g(d, m)$
means that, for fixed positive integer $d$,
$$f(d, m)\le c_d\cdot g(d,m)$$
holds for all positive integers $m$, where $c_d$ is a suitable constant depending only on $d$. In
1980, Arnold \cite{arno80} studied the values of $v(2,m)$ and proved
$$m^{1\over 3}\ll \log v(2,m)\ll m^{1\over 3}\log m.\eqno (1)$$
In 1992, B\'ar\'any and Pach \cite{bara92} improved Arnold's
upper bound to
$$\log v(2,m)\ll m^{1\over 3},\eqno(2)$$
B\'ar\'any and Vershik \cite{bara92'} generalized (2) to
$d$ dimensions by proving
$$\log v(d,m)\ll m^{{d-1}\over {d+1}}.\eqno (3)$$
Although Arnold at the end of his paper \cite{arno80} where he proved (1) wrote
\lq\lq {\it In $\mathbb{Z}^d$, $1/3$ is probably replaced by $(d-1)/(d+1)$. Proof of the
lower bound: let $x_1^2+\ldots +x_{d-1}^2\le x_d\le A$.}"
the problem whether 
$$\log v(d,m)\gg m^{{d-1}\over {d+1}}$$
is still open.

Let $v^\ast(d,m)$ denote the number of different classes of the $d$-dimensional
centrally symmetric convex lattice polytopes $P$ with $v(P)=m/d!$, let $\kappa(d,w)$
denote the number of different classes of $d$-dimensional convex lattice polytopes
$P$ with $|P|=w$, let $\kappa^\ast(d,w)$ denote the number of different classes of
$d$-dimensional centrally symmetric convex lattice polytopes $P$ with $|P|=w$,
and let $\kappa '(d,w)$ denote the number of different classes of $d$-dimensional
convex lattice polytopes $P$ with $|P|=w$ and ${\rm int}(P)\cap \mathbb{Z}^d\not=\emptyset $.
Then we have $v^\ast(d,m)=0$ whenever $m$ is odd and $\kappa^\ast(d,w)=0$ if $w$ is even.
Therefore in this paper we assume that the $m$ in $v^\ast(d,m)$ is even and the $w$ in
$\kappa^\ast(d,w)$ is odd.

\medskip
\noindent
{\bf Remark 1.} As usual, in this paper centrally symmetric convex lattice polytopes are those centered
at lattice points. In this sense, the unit cube $\{ {\bf x}\in \mathbb{E}^d:\ 0\le x_i\le 1\}$ is not a
centrally symmetric convex lattice polytope, though it is a convex lattice polytope and is centrally symmetric.

\medskip
Recently, Liu and Zong \cite{zong11} studies Arnold's problem for the centrally symmetric
lattice polygons and the classification problem for convex lattice polytopes of given
cardinality by proving
$$m^{1\over 3}\ll \log v^\ast(2, m)\ll m^{1\over 3},$$
$$w^{1\over 3}\ll \log \kappa(2, w)\ll w^{1\over 3},$$
$$w^{1\over 3}\ll \log \kappa^\ast(2, w)\ll w^{1\over 3},$$
$$\kappa(d,w)=\infty , \quad if\ w\ge d+1\ge 4,$$
$$\log \kappa '(d,w)\ll w^{{d-1}\over {d+1}}\eqno (4)$$
and
$$\log \kappa^\ast (d, w)\ll w^{{d-1}\over {d+1}}.\eqno (5)$$

In Section 2 of this paper we introduce and study unimodular groups acting on convex
lattice polytopes. In particular, the orders of these groups are estimated. In Section
3, by applying the results obtained in Section 2, we prove the following result:

\medskip\noindent
{\bf Theorem 1.} {\it Let $\kappa^\ast(d,w)$ denote the number of different classes of
$d$-dimensional centrally symmetric convex lattice polytopes $P$ with $|P|=w$, then}
$$\log \kappa^* (d,w) \gg w^{{d-1}\over {d+1}}.$$

\smallskip
This theorem, together with (4) and (5), produces the following consequences (when $w$ is even, to deduce
Theorem 3 needs a little extra care since it is no longer symmetric):

\medskip
\noindent
{\bf Theorem 2.} {\it Let $\kappa^\ast(d,w)$ denote the number of different classes of
$d$-dimensional centrally symmetric convex lattice polytopes $P$ with $|P|=w$, then}
$$ w^{{d-1}\over {d+1}}\ll \log\kappa^* (d,w)\ll w^{{d-1}\over {d+1}}.$$

\medskip
\noindent
{\bf Theorem 3.} {\it Let $\kappa '(d,w)$ denote the number of different classes of
$d$-dimensional convex lattice polytopes $P$ with $|P|=w$ and ${\rm int}(P)\cap
\mathbb{Z}^d\not=\emptyset$, then}
$$ w^{{d-1}\over {d+1}}\ll \log\kappa ' (d,w)\ll w^{{d-1}\over {d+1}}.$$

\section*{2. Unimodular groups of convex lattice polytopes}

In this section we introduce and study unimodular groups acting on convex lattice polytopes. In particular, Lemma 3 will be essential for our proof of Theorem 1.

As usual, a unimodular transformation $\sigma ({\bf x})$ of $\mathbb{E}^d$ is a $\mathbb{Z}^d$-preserving
affinely linear transformation, i.e.,
$$\sigma ({\bf x})={\bf x}U+{\bf v},$$
where $U$ is a $d\times d$ integral matrix satisfying $|det (U)|=1$ and ${\bf v}$ is an integral vector.
In particular, if $U$ also satisfies $UU'=I$, where $U'$ is the transpose of $U$ and $I$ is the
$d\times d$ unit matrix, we call $\sigma ({\bf x})$ an orthogonal unimodular transformation. It is
known in linear algebra that an orthogonal unimodular keeps the Euclidean distances unchanged.

Let $\sigma_1$ and $\sigma_2$ be two unimodular transformations in $\mathbb{E}^d$. It is known in linear
algebra that both $\sigma_1\cdot \sigma_2$ and $\sigma_1^{-1}$ are unimodular
transformations. Therefore, all unimodular transformations in $\mathbb{E}^d$ form a multiplicative group.
We denote it by $\mathbb{G}_d$. Similarly, all orthogonal unimodular transformations in $\mathbb{E}^d$
form a subgroup of $\mathbb{G}_d$. We denote it by $\mathbb{G}'_d$.

The group $\mathbb{G}_d$ is different from ${\rm GL}(d,\mathbb{Z})$. But, they are closely related. In fact, we have
$$\mathbb{G}_d \cong \left\{ \left(\begin{array}{cc}
A& O\\
V& 1
\end{array}
\right):\ A\in {\rm GL} (d, \mathbb{Z}),\ V\in \mathbb{Z}^d,\ O=(0,0,\ldots , 0)'\right\}
\subseteq  {\rm GL}(d+1, \mathbb{Z}).$$

Let $P$ be a convex lattice polytope in $\mathbb{E}^d$, and let $\mathbb{P}_d$ denote the family of all
$d$-dimensional convex lattice polytopes. We define
$$\sigma (P)=\{ \sigma ({\bf x}):\ {\bf x}\in P\}.$$
Clearly, $\sigma (P)$ is a convex lattice polytope as well. Then we define
$$G(P)=\{ \sigma \in \mathbb{G}_d:\ \sigma (P)=P\}$$
and
$$G'(P)=\{ \sigma \in \mathbb{G}'_d:\ \sigma (P)=P\}.$$
Both $G(P)$ and $G'(P)$ are finite subgroups of $\mathbb{G}_d$, and $G'(P)$
is a subgroup of $G(P)$. We call $G(P)$ the {\it unimodular group} of $P$ and call $G'(P)$ the
{\it orthogonal unimodular group} of $P$.

It is easy to see that both $\mathbb{G}_d$ and $\mathbb{G}'_d$ act on $\mathbb{P}_d$, $G(P)$ is the
stabilizer of $P$ in $\mathbb{G}_d$, and $G'(P)$ is the stabilizer of $P$ in $\mathbb{G}'_d$. Therefore,
we have

\medskip\noindent
{\bf Lemma 1.} {\it If $P\in \mathbb{P}_d$ and $\sigma \in \mathbb{G}_d$, then we have
$$G(\sigma (P))=\sigma G(P)\sigma^{-1} .$$
If $\tau \in \mathbb{G}'_d$, then we have}
$$G'(\tau (P))=\tau G'(P)\tau^{-1} .$$

\noindent
{\bf Lemma 2.} {\it Let $O_d$ denote the multiplicative group of orthogonal unimodular
transformations of $\mathbb{E}^d$ which keep the origin fixed. Then, we have}
$$|O_d|=2^d\cdot d!.$$

\medskip
If $\sigma \in O_d$, then we have
$$\sigma ({\bf e}_i)\in \{ \pm {\bf e}_1, \pm {\bf e}_2, \ldots , \pm {\bf e}_d\},$$
by which one can easily deduce the lemma. In fact, $O_d$ is the multiplicative group of the
$d\times d$ orthogonal integral matrices.

\bigskip\noindent
{\bf Corollary 1.} {\it For any $d$-dimensional centrally symmetric convex lattice polytope $P$,
$\left| G'(P)\right|$ is a divisor of $2^d\cdot d!.$}

\medskip
We have two basic problems about the unimodular groups of convex lattice polytopes.

\medskip\noindent
{\bf Problem 1.} {\it Determine the values of}
$$\max_{P\in \mathbb{P}_d} \ \{ |G(P)|\}.$$

\medskip\noindent
{\bf Remark 2.} It was proved by Minkowski \cite{mink87} that, there is a constant $c_d$ depends only on $d$,
$$|G|\le c_d$$
holds for any finite subgroup $G$ of ${\rm GL}(d,\mathbb{Z})$. According to \cite{gura06} and \cite{kuzm02}, in an unpublished manuscript W. Feit proved (based on an unfinished manuscript of B. Weisfeiler) that, if $d> 10$
and $G$ is a finite subgroup of ${\rm GL}(d,\mathbb{Z})$, then
$$|G|\le 2^d\cdot d!.$$
Feit's result implies that, when $d>10$,
$$|G(C)|\le 2^d\cdot d!$$
holds for all $d$-dimensional centrally symmetric lattice polytopes $C$.

\medskip\noindent
{\bf Problem 2.} {\it Is it true that, for any finite subgroup $G$ of $\mathbb{G}_d$, there is a $P\in
\mathbb{P}_d$ such that $G(P)\cong G$? Is it true that, for any finite subgroup $G$ of ${\rm GL}(d,\mathbb{Z})$, there is a $d$-dimensional centrally symmetric convex lattice polytope $P$ satisfying $G(P)\cong G$?}

\medskip
These problems show the close relation between $d$-dimensional convex lattice polytopes and the finite subgroups of ${\rm GL}(d,\mathbb{Z})$ and ${\rm GL}(d+1,\mathbb{Z})$.

\medskip
Let $m$ be a positive integer and let $\rho$ be a real number satisfying $1\le \rho \le \infty$.
We define
$$P_{d,m,\rho }={\rm conv}\left\{ {\bf z}\in \mathbb{Z}^d:\ \left(\sum_{i=1}^d|z_i|^\rho\right)^{1/\rho } \le m \right\},$$
where $\left(\sum |z_i|^\infty \right)^{1/\infty }=\max \{ |z_i|\}.$ One can easily verify that $P_{d,m,\rho }$ is a $d$-dimensional centrally symmetric convex lattice polytope.
In particular, $P_{d,m,1}$ is a lattice cross-polytope and $P_{d,m,\infty }$ is a lattice cube.

\medskip\noindent
{\bf Lemma 3.} {\it When $d$ and $m$ are positive integers and $\rho$ is a positive number satisfying
$1\le \rho \le \infty$, we have
$$G\left(P_{d,m,\rho }\right)=G'\left(P_{d,m,\rho }\right)=O_d$$
and
$$\left|G\left(P_{d,m,\rho }\right)\right|=\left|G'\left(P_{d,m,\rho }\right)\right|=2^d\cdot d!.$$}

\noindent
{\bf Proof.} First of all, since $P_{d,m,\rho }$ is centrally symmetric and centered at the origin,
we have
$$\sigma ({\bf o})={\bf o}\eqno (6)$$
for all $\sigma \in G(P_{d,m,\rho }).$

Let ${\bf v}$ be a primitive integral vector in $\mathbb{Z}^d$ and let $P$ be a centrally symmetric
convex lattice polytope in $\mathbb{E}^d$. We define
$$L(P,{\bf v})=\{ z{\bf v}:\ z\in \mathbb{Z}\}\cap P$$
and
$$\ell (P)=\max_{\bf v}\{ |L(P, {\bf v})|\},$$
where the maximum is over all primitive integral vectors in $\mathbb{Z}^d$. Recall that $\{ {\bf e}_1,
{\bf e}_2, \ldots , {\bf e}_d\}$ is an orthonormal basis of $\mathbb{E}^d$. We consider two cases
as follows:

\medskip\noindent
{\bf Case 1.} $\rho < \infty $. Notice that
$$\sum_{i=1}^d|mv_i|^\rho =m^\rho \sum_{i=1}^d|v_i|^\rho ,$$
it can be easily deduce that
$$\ell (P_{d,m,\rho })=2m+1$$
and
$$\left| L(P_{d,m,\rho }, {\bf v})\right|=2m+1$$
holds if and only if ${\bf v}=\pm {\bf e}_i$ for some index $i$. Thus, for any $\sigma \in G(P_{d,m,\rho })$,
we have
$$\left\{\sigma ({\bf e}_1), \sigma ({\bf e}_2), \ldots ,\sigma ({\bf e}_d)\right\}\subset
\left\{ \pm {\bf e}_1, \pm {\bf e}_2, \ldots , \pm {\bf e}_d\right\}.\eqno (7)$$

\medskip\noindent
{\bf Case 2.} $\rho = \infty $. In this case $P_{d,m,\infty }$ is a $d$-dimensional cube. It has $2d$ facets $\pm F_1$, $\pm F_2$,
$\ldots $, $\pm F_d$, each is a $(d-1)$-dimensional cube. The centers of the facets are
$\pm m{\bf e}_1$, $\pm m{\bf e}_2$, $\ldots $, $\pm m{\bf e}_d$. If $\sigma \in G(P_{d,m,\infty })$,
we have
$$\left\{\sigma (F_1), \sigma (F_2),\ldots, \sigma (F_d)\right\}\subset \{ \pm F_1, \pm F_2, \ldots ,
\pm F_d\},$$
$$\left\{\sigma (m{\bf e}_1), \sigma (m{\bf e}_2),\ldots, \sigma (m{\bf e}_d)\right\}\subset
\{ \pm m{\bf e}_1, \pm m {\bf e}_2, \ldots , \pm m{\bf e}_d\}$$
and therefore
$$\left\{\sigma ({\bf e}_1), \sigma ({\bf e}_2), \ldots , \sigma ({\bf e}_d)\right\}\subset
\{ \pm {\bf e}_1, \pm {\bf e}_2, \ldots , \pm {\bf e}_d\}.\eqno(8)$$

Assume that the unimodular transformation $\sigma $ is defined by
$$\sigma ({\bf x})={\bf x}U+{\bf b}.$$
It follows by (6) that ${\bf b}={\bf o}$. In both cases, since $U$ is nonsingular, by (7) and (8) we get
$$\sum_{j=1}^d|u_{ij}|=1,\quad i=1, 2, \ldots , d$$
and
$$\sum_{i=1}^d|u_{ij}|=1,\quad j=1, 2, \ldots , d.$$
Thus, we obtain
$$G(P_{d,m,\rho })\subseteq O_d.\eqno (9)$$

On the other hand, it is easy to verify that
$$O_d\subseteq G'(P_{d,m,\rho }).\eqno (10)$$
As a conclusion of (9) and (10) we get
$$O_d\subseteq G'(P_{d,m,\rho })\subseteq G(P_{d,m,\rho })\subseteq O_d$$
and finally
$$G(P_{d,m,\rho })=G'(P_{d,m,\rho })=O_d.$$

The second assertion of the lemma follows from Lemma 2. The lemma is proved. \hfill{$\square $}

\medskip\noindent
{\bf Remark 3.} Let $S_d$ denote the $d$-dimensional lattice simplex with vertices ${\bf e}_1$,
${\bf e}_2$, $\ldots$, ${\bf e}_d$ and ${\bf o}$. Then we have $|G'(S_d)|=d!$ and $|G(S_d)|=(d+1)!.$ This example
shows that $|G'(P)|$ and $|G(P)|$ can be different.

\section*{3. Proof of Theorem 1}

In this section we study the classification problem for convex lattice polytopes of given cardinality.
In particular, Theorem 1 will be proved.

First, as illustrated by Figure 1, for a positive integer $r$ we define
\begin{eqnarray*}
K_{d,r}&=&\left\{ {\bf x}\in \mathbb{E}^d:\ x_d\ge 0,\ x_d+\sum_{i=1}^{d-1}x_i^2\le r^2\right\},\\
B_{d,r}&=&\left\{ {\bf x}\in \mathbb{E}^d:\ x_d=0,\ \sum_{i=1}^{d-1}x_i^2\le r^2\right\}
\end{eqnarray*}
and
$$P_{d,r}= {\rm conv}\left\{K_{d,r}\cap \mathbb{Z}^d\right\}.$$

$$\psannotate{\psboxto(0cm;6.3cm){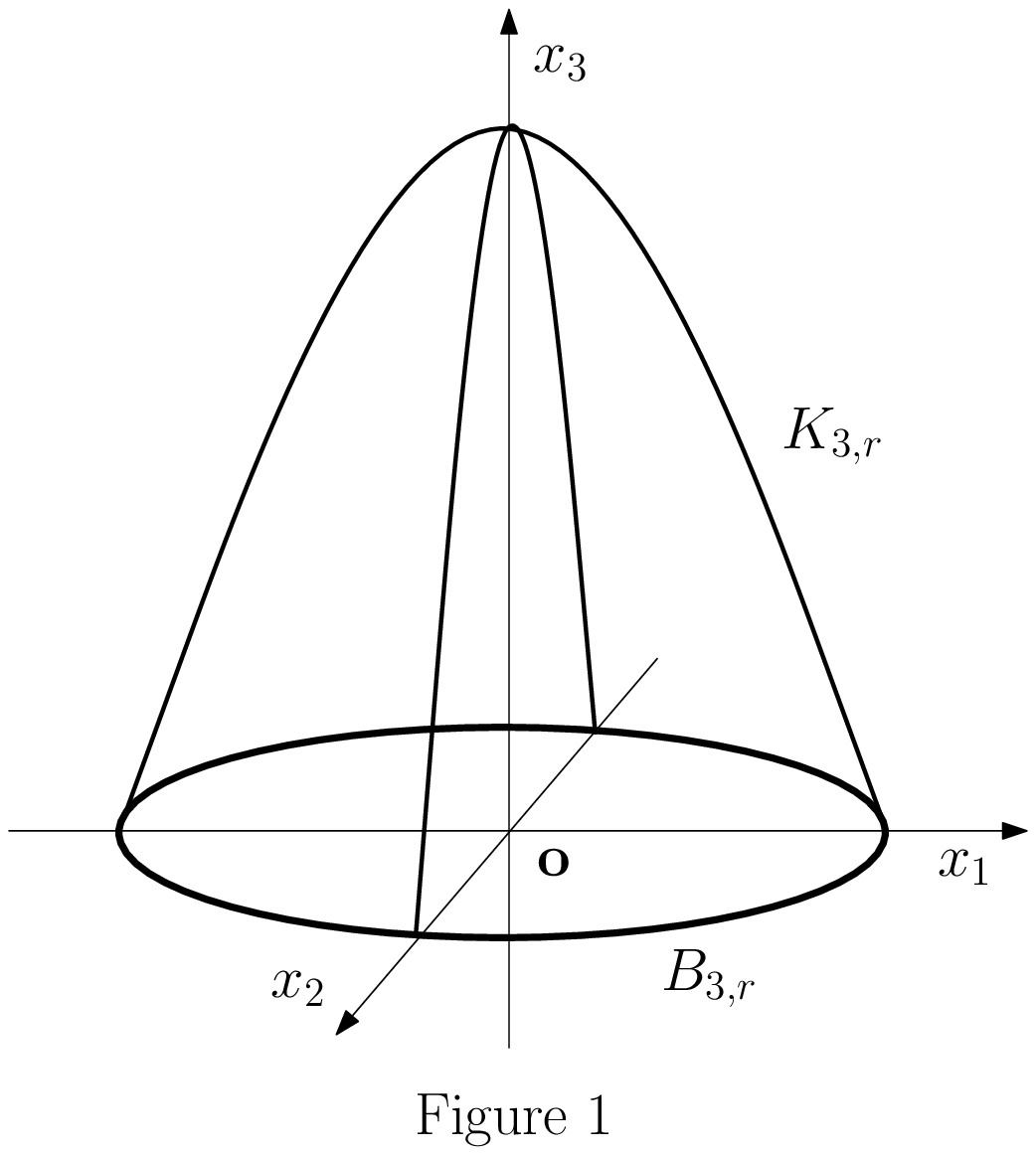}}{}$$

It is easy to compute that
$$v(K_{d,r})=\int_0^{r^2}{{\pi^{{d-1}\over 2}}\over {\Gamma ({{d+1}\over 2})}}\cdot (r^2-x)^{{d-1}\over 2}dx
=c_1(d)\cdot r^{d+1}\eqno (11)$$
and
$$s'(K_{d,r})=\int_0^{r^2}{{(d-1)\pi^{{d-1}\over 2}}\over {\Gamma ({{d+1}\over 2})}}\cdot (r^2-x)^{{d-2}\over 2}dx
=(d-1)c_1(d)\cdot r^d,\eqno (12)$$
where $s'(K_{d,r})$ denotes the surface area of $K_{d,r}$ without the base and $c_1(d)$ is a constant that depends only on $d$.

Next, we define
$$C_{d,r}^1=\left\{ {\bf x}\in \mathbb{E}^d:\ 0\le x_d\le r^2,\ \sum_{i=1}^{d-1}x_i^2\le r^2\right\},$$
$$C_{d,r}^2=\left\{ {\bf x}\in \mathbb{E}^d:\ -1\le x_1\le 1,\ 0\le x_d\le r^2,\ x_d+\sum_{i=2}^{d-1}x_i^2\le r^2\right\}$$
and their intersection
$$C_{d,r}=C_{d,r}^1\cap C_{d,r}^2.$$

$$\psannotate{\psboxto(0cm;6.2cm){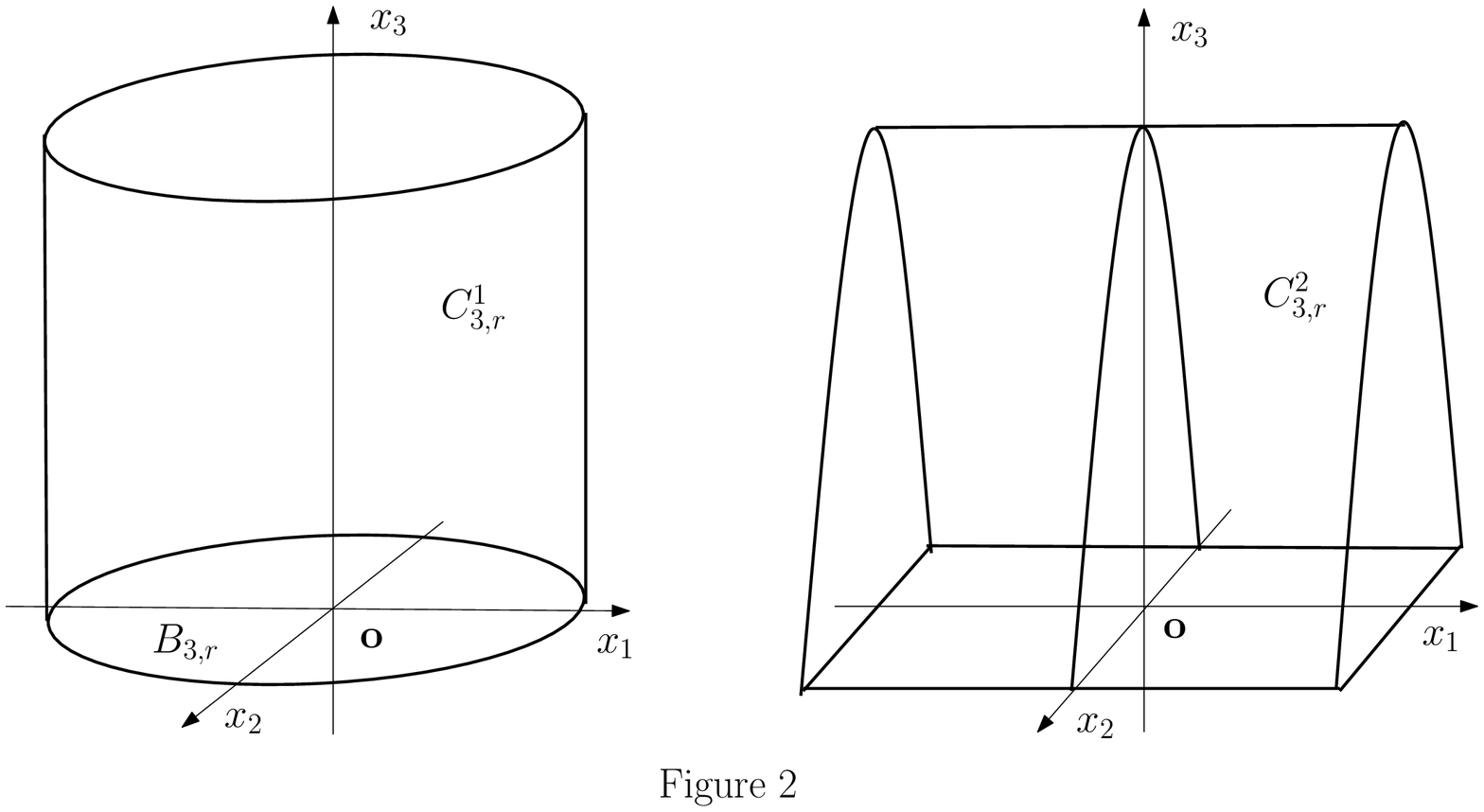}}{}$$

In fact, $C^1_{d,r}$ is a cylinder over a base $B_{d,r}$ and $C^2_{d,r}$ is a cylinder over a base
$$\left\{ {\bf x}\in\mathbb{E}^d:\ x_d\ge 0, \ x_1=0,\ \sum_{i=2}^{d-1}x_i^2+x_d\le r^2\right\},$$
as shown in Figure 2.

$$\psannotate{\psboxto(0cm;6.3cm){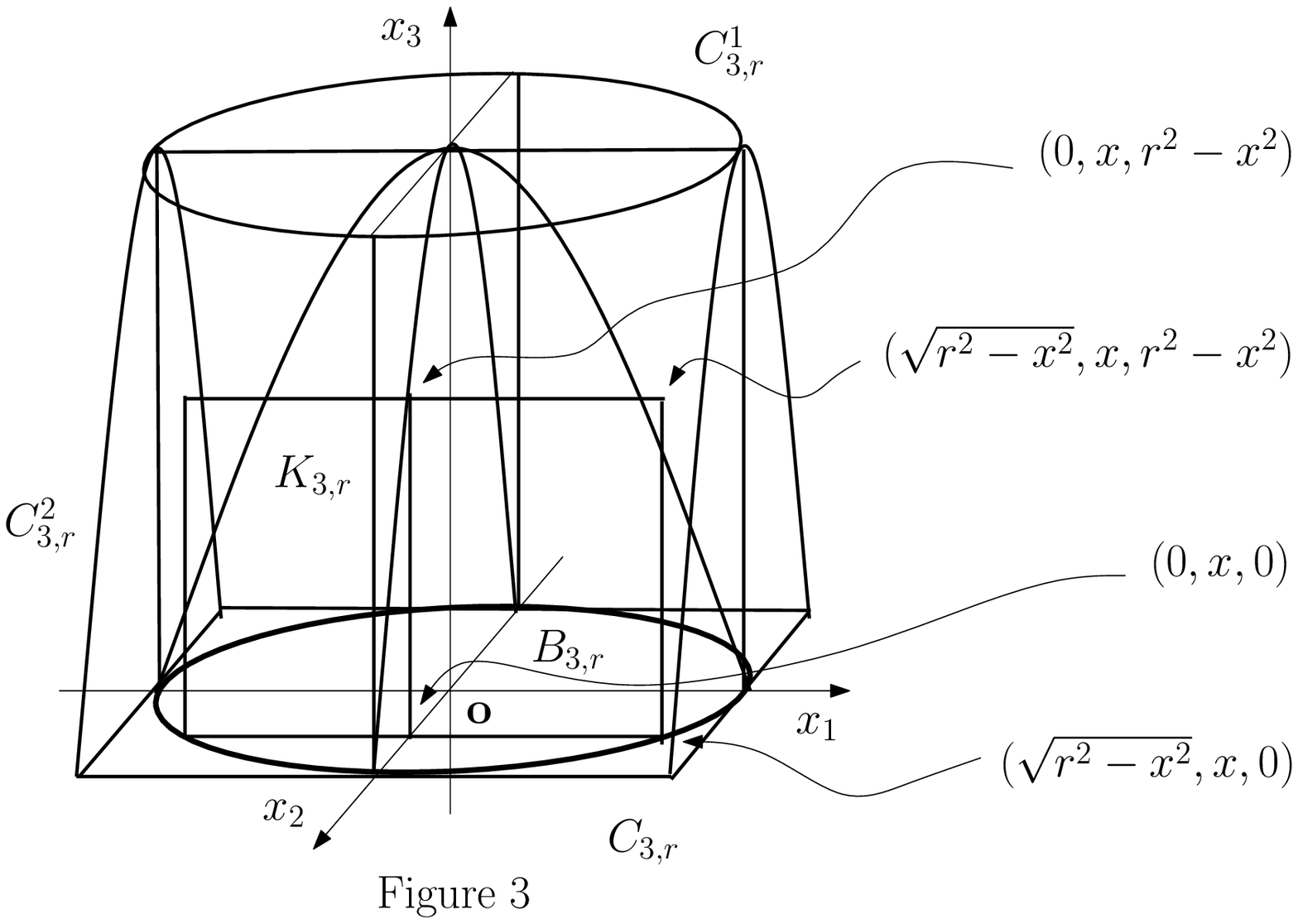}}{}$$

Notice that (as illustrated in Figure 3), if $$x^2=\sum_{i=2}^{d-1}x_i^2\eqno (13)$$ and
$$(0, x_2, x_3, \ldots , x_{d-1}, 0)\in C_{d,r},\eqno (14)$$ then we have
$$(x_1, x_2, \ldots , x_{d-1}, x_d)\in C_{d,r}$$
provided $|x_1|\le \sqrt{r^2-x^2}$ and $0\le x_d\le r^2-x^2.$ In fact, all the points
satisfying both (13) and (14) together form a $(d-2)$-dimensional sphere of radius $x$ which has area measure
$${{(d-2)\pi^{{d-2}\over 2}}\over {\Gamma ({d\over 2})}}\cdot x^{d-3}.$$
Thus, we get
\begin{eqnarray*}
\hspace{1cm}v(C_{d,r})&=&\int_0^r{{2(d-2)\pi^{{d-2}\over 2}}\over {\Gamma ({d\over 2})}}\cdot x^{d-3}(r^2-x^2)^{3\over 2}dx\\
&=&{{2(d-2)\pi^{{d-2}\over 2}}\over {\Gamma ({d\over 2})}}\cdot r^{d+1}\int_0^{\pi\over 2}\sin^{d-3}\theta
\cos^4\theta d\theta\\
&=&{{2(d-2)\pi^{{d-2}\over 2}}\over {\Gamma ({d\over 2})}}\cdot r^{d+1}\int_0^{\pi\over 2}\left(
\sin^{d-3}\theta -2\sin^{d-1}\theta +\sin^{d+1}\theta\right)d\theta\\
&=&c_2(d)\cdot r^{d+1} \hspace{7.8cm}(15)
\end{eqnarray*}
and
$$s'(C_{d,r})\le s'(C_{d,r}^1)+s'(C_{d,r}^2)\le c_3(d)\cdot r^d, \eqno (16)$$
where $c_2(d)$ and $c_3(d)$ are constants depending only on $d$.
It follows by $K_{d,r}\subset C_{d,r}$ that
$$c_1(d)<c_2(d).\eqno (17)$$

For convenience, we write $W=\{ (x_1, x_2, \ldots , x_d):\ |x_i|\le {1\over 2}\}$ and $B=\{ (x_1, x_2, \ldots , x_d):\ \sum x_i^2\le 1\}$. By convexity we have
$$v\left( C_{d,r}\right)-\mbox{$1\over 2$}\sqrt{d}\cdot  s'(C_{d,r})\le v\left(W+C_{d,r}\cap \mathbb{Z}^d\right)\le v\left(C_{d,r}+ \mbox{$1\over 2$}\sqrt{d}\ B\right).$$
Then, by (15) and (16) we get
$$v\left(C_{d,r}\right)-\mbox{$1\over 2$}\sqrt{d}\cdot c_3(d)\cdot r^d\le \left|C_{d,r}\cap \mathbb{Z}^d\right|\le v\left( C_{d,r}\right)+
O (r^d)$$
and therefore
$$\left|C_{d,r}\cap \mathbb{Z}^d\right|\sim v\left( C_{d,r}\right).\eqno (18)$$
Similarly, by (11) and (12) we get
$$\left|P_{d,r}\right|\sim v\left( K_{d,r}\right).\eqno (19)$$

Next, we define
$$Q_{d,r}=\left( {\rm int}\left(C^1_{d,r}\right)\cap C_{d,r}^2\right)\cup B_{d,r},\eqno(20)$$
$$H_{d,r}={\rm conv}\left\{ {\bf z}\in Q_{d,r}\cap \mathbb{Z}^d: z_1\le 0\right\}\eqno (21)$$
and
$$H'_{d,r}={\rm conv}\left\{ {\bf z}\in K_{d,r}\cap \mathbb{Z}^d: z_1\le 0\right\}.\eqno (22)$$
By (18) and (19), one can deduce
$$|H_{d,r}|\sim {1\over 2}\cdot \left|C_{d,r}\cap \mathbb{Z}^d\right|\sim {1\over 2}\cdot v(C_{d,r})={{c_2(d)}\over 2} \cdot r^{d+1}\eqno (23)$$
and
$$|H'_{d,r}|\sim {1\over 2}\cdot \left|P_{d,r} \right| \sim {1\over 2}\cdot v(K_{d,r})={{c_1(d)}\over 2} \cdot r^{d+1}.\eqno (24)$$

\smallskip\noindent
{\bf Remark 4.} Let $L({\bf x})$ denote the line defined by $\{ {\bf x}+\lambda {\bf e}_d:\ \lambda
\in \mathbb{R}\}$. When ${\bf z}$ is a lattice point on the boundary of $B_{d,r}$, we have
$$\left| L({\bf z})\cap Q_{d,r}\cap \mathbb{Z}^d\right| =1.$$

\medskip
Now, we introduce a technical lemma which is useful in the proof of Theorem 1.

\medskip\noindent
{\bf Lemma 4.} {\it When $r$ is a sufficiently large integer, for any integer $k$ satisfying
$0\le k\le 3r |B_{d,r+1}\cap \mathbb{Z}^d|$, there is a convex lattice polytope $P$ that satisfies
$$H'_{d,r}\subseteq P\subseteq H_{d,r}$$
and}
$$|P|=|H'_{d,r}|+k.$$

\medskip\noindent
{\bf Proof.} It is well-known that
\begin{eqnarray*}
\left| B_{d,r+1}\cap \mathbb{Z}^d\right| &=& {{\pi^{{d-1}\over 2}}\over {\Gamma ({{d+1}\over 2})}}\cdot (r+1)^{d-1}
+O \left( (r+1)^{d-2}\right).\\
&=& {{\pi^{{d-1}\over 2}}\over {\Gamma ({{d+1}\over 2})}}\cdot r^{d-1}
+O \left( r^{d-2}\right).
\end{eqnarray*}
By (23), (24) and (17), when $r$ is sufficiently large, we get
\begin{eqnarray*}
\hspace{2.5cm}|H_{d,r}|-|H'_{d,r}|&\ge & {1\over 4}\cdot (c_2(d)-c_1(d))\cdot r^{d+1}\\
&\ge &c_4(d)\cdot r^2\cdot |B_{d,r+1}\cap \mathbb{Z}^d|\\
&\ge & 3r |B_{d,r+1}\cap\mathbb{Z}^d|,\hspace{4.4cm}(25)
\end{eqnarray*}
where $c_4(d)$ is a constant that depends only on $d$.

For convenience, we write $P_0=H_{d,r}$ and let $\overline{P}$ denote the set of the vertices of $P$.
If ${\bf v}_0\in \overline{P_0}\setminus H'_{d,r}$, we define
$$P_1={\rm conv}\left\{ (P_0\cap \mathbb{Z}^d)\setminus \{ {\bf v}_0\}\right\}.$$
Inductively, if $P_i$ has been defined and ${\bf v}_i\in \overline{P_i}\setminus H'_{d,r}$, we construct
$$P_{i+1}={\rm conv} \left\{ (P_i\cap \mathbb{Z}^d)\setminus \{ {\bf v}_i\}\right\}.$$
Thus, we have constructed a finite sequence of convex lattice polytopes $P_0$, $P_1$, $P_2$,
$\ldots $, $P_\ell =H'_{d,r}$ which satisfies both
$$P_0\supset P_1\supset \ldots \supset P_{\ell -1}\supset P_\ell =H'_{d,r}$$
and
$$|P_i|-|P_{i+1}|=1,\quad i=0, 1, 2, \ldots , \ell -1.$$
By (25) it follows that
$$\ell \ge 3r\left|B_{d,r+1}\cap \mathbb{Z}^d\right|.$$
The assertion is proved. \hfill{$\square $}

\bigskip\noindent
{\bf Proof of Theorem 1.} First, we recall that
$$P_{d,r}={\rm conv}\left\{K_{d,r}\cap \mathbb{Z}^d\right\}.\eqno(26)$$
Let $w$ be a large odd integer and let $r$ be the integer satisfying
$$|P_{d,r}|\le {w\over 2}<|P_{d,r+1}|.\eqno (27)$$
By (11) and (19) we get
$$r^{d+1}\ll w\ll r^{d+1}.\eqno (28)$$

We write
$$P_{d,r+1}^1=\left\{ {\bf x}\in P_{d,r+1}:\ x_d\ge 2r+1\right\}$$
and
$$P_{d,r+1}^2=\left\{ {\bf x}\in P_{d,r+1}:\ x_d\le 2r\right\}.$$
It is easy to see that
$$|P_{d,r}|=\left|P^1_{d,r+1}\right|$$
and therefore
$$|P_{d,r+1}|-|P_{d,r}|=\left|P^2_{d,r+1}\right|\le (2r+1)\left|B_{d,r+1}\cap \mathbb{Z}^d\right|.\eqno (29)$$

We write
$$u=w-2|P_{d,r}|+\left|B_{d,r}\cap \mathbb{Z}^d\right|.\eqno (30)$$
By (27) and (29) we get
$$u<2\left( |P_{d,r+1}|-|P_{d,r}|\right)+\left|B_{d,r}\cap \mathbb{Z}^d\right|
\le 5r \left|B_{d,r+1}\cap \mathbb{Z}^d\right|. \eqno(31)$$

Let $V'_{d,r}$ denote the set of the vertices ${\bf v}$ of $P_{d,r}$ satisfying both $v_d\not= 0$ and
$v_1\ge 1$, and let $L({\bf x})$ denote the line $\{ {\bf x}+\lambda {\bf e}_d:\ \lambda \in \mathbb{R}\}$
as defined in Remark 4. By convexity, for all ${\bf z}\in B_{d,r}\cap \mathbb{Z}^d$ with $z_1\ge 1$ we have
$$\left| L({\bf z})\cap V'_{d,r}\right|\le 1.$$
Thus we get
$$\left|V'_{d,r}\right|< {1\over 2} \left|B_{d,r}\cap \mathbb{Z}^d\right|\eqno (32)$$
and
\begin{eqnarray*}
\hspace{2.7cm}\left|V'_{d,r}\right|&\ge &{1\over 2}\cdot \left( \left|{\rm int}(B_{d,r})\cap \mathbb{Z}^d\right|-\left|B_{d-1,r}\cap \mathbb{Z}^d\right|\right)\\
&\ge & {1\over 3}\cdot \left( {{\pi^{{d-1}\over 2}}\over {\Gamma ({{d+1}\over 2})}}\cdot r^{d-1}
-{{\pi^{{d-2}\over 2}}\over {\Gamma ({d\over 2})}}\cdot r^{d-2}\right)\\
&\gg &r^{d-1}.\hspace{7cm}(33)
\end{eqnarray*}

With these preparations, we proceed to construct the expected convex lattice polytopes.

\medskip\noindent
{\bf Step 1.} Let ${\bf v}_1$, ${\bf v}_2$, $\ldots $, ${\bf v}_j$ be $j$ points in $V'_{d,r}$ and
define
$$P'_{d,r}={\rm conv}\left\{ P_{d,r}\setminus \{ {\bf v}_1, {\bf v}_2, \ldots , {\bf v}_j\}\right\}.\eqno (34)$$
We have
$$\left|P'_{d,r}\right|=|P_{d,r}|-j.\eqno (35)$$
By (31) and (32) we get
\begin{eqnarray*}
{u\over 2}+j&\le &{5\over 2}r\left| B_{d,r+1}\cap \mathbb{Z}^d\right|+{1\over 2}\left| B_{d,r}
\cap \mathbb{Z}^d\right|\\
&<&3r \left|B_{d,r+1}\cap \mathbb{Z}^d\right|.
\end{eqnarray*}
According to Lemma 4, there is a convex lattice polytope $P$ satisfies both
$$H'_{d,r}\subseteq P\subseteq H_{d,r}\eqno (36)$$
and
$$|P|=\left|H'_{d,r}\right|+{u\over 2}+j.\eqno (37)$$

\medskip\noindent
{\bf Step 2.} We construct
$$P'=P\cup P'_{d,r}.\eqno (38)$$
Let $i$ be an integer and let $H_i$ denote the hyperplane $\{ (x_1, x_2, \ldots , x_d):\ x_1=i\}$. By (34) and (36)
it is easy to see that $P'\cap H_i-i{\bf e}_i$ is a subset of $P'\cap H_0$. Therefore the convexity of $P'$ at both sides of $H_0$ guarantees the convexity of $P'$. In other words, we have
$$P'={\rm conv}\{ P\cup P'_{d,r}\}.$$
By (35) and (37) we get
$$|P'|=|P_{d,r}|+{u\over 2}.\eqno (39)$$

\medskip\noindent
{\bf Step 3.} We define
$$P_{{\bf v}_1, \ldots, {\bf v}_j}=P'\cup \{ -P'\}.\eqno (40)$$
Clearly $P_{{\bf v}_1, \ldots, {\bf v}_j}$ is a centrally symmetric convex lattice polytope
centered at the origin. By (39) and (30) we get
\begin{eqnarray*}
\left|P_{{\bf v}_1, \ldots, {\bf v}_j}\right|&=&2\left( |P_{d,r}|+{u\over 2}\right)-\left|B_{d,r}\cap \mathbb{Z}^d\right|\\
&=&2|P_{d,r}|+u-\left| B_{d,r}\cap \mathbb{Z}^d\right|\\
&=&2|P_{d,r}|+w-2|P_{d,r}|+\left| B_{d,r}\cap \mathbb{Z}^d\right|-\left| B_{d,r}\cap \mathbb{Z}^d\right|\\
&=& w.\end{eqnarray*}

\medskip\noindent
{\bf Step 4.} Taking all possible subsets of $V'_{d,r}$, we get $2^{|V'_{d,r}|}$ centrally symmetric convex
lattice polytopes of cardinality $w$. For convenience, we enumerate them by $P_1$, $P_2$, $\ldots $, $P_{2^{|V'_{d,r}|}}$ and denote the whole family by $\mathcal{F}$.

\medskip
Now, we study the equivalence relation among $\mathcal{F}$.

Recalling the definitions of $L(P,{\bf v})$ and $\ell (P)$ in the proof of Lemma 3, for any
$P_i\in \mathcal{F}$ we have
$$\ell (P_i)=2r^2+1$$
and
$$L(P_i,{\bf v})=2r^2+1$$
holds if and only if ${\bf v}=\pm {\bf e}_d$. Therefore, if $\sigma (P_i)=P_j$ holds for some unimodular transformation $\sigma $, we have $\sigma ({\bf o})={\bf o}$ and $\sigma ({\bf e}_d)\in \{{\bf e}_d, -{\bf e}_d\}.$

Let $H$ be a $(d-1)$-dimensional hyperplane which contains the origin of $\mathbb{E}^d$, but not ${\bf e}_d$.
By projecting $P_i\cap H\cap \mathbb{Z}^d$ onto the plane $\{ {\bf x}\in \mathbb{E}^d:\ x_d=0\}$, keeping
(20), (21), (26), (34), (36), (38), (40) and Remark 4 in mind, it follows that
$$\left|P_i\cap H\cap \mathbb{Z}^d\right| \le\left|B_{d,r}\cap \mathbb{Z}^d\right|,$$
where the equality holds if and only if
$$H=\left\{ {\bf x}\in \mathbb{E}^d:\ x_d=0\right\}.$$
Thus, we get
$$\sigma \left(B_{d,r}\cap \mathbb{Z}^d\right)=B_{d,r}\cap \mathbb{Z}^d$$
and therefore
$$\sigma\in G\left({\rm conv}\left\{ B_{d,r}\cap \mathbb{Z}^d\right\}\right)=G(P_{d-1, r,2}).$$
Consequently, by the $\rho =2$ case of Lemma 3, we get
$$\kappa ^*(d,w)\ge {{|\mathcal{F}|}\over {2|G(P_{d-1,r,2})|}}={{2^{|V'_{d,r}|}}\over
{2^d\cdot (d-1)!}}.$$
By (33) and (28), we deduce
$$\log \kappa ^*(d,w)\gg |V'_{d,r}|\gg r^{d-1}\gg w^{{d-1}\over {d+1}}.$$
The proof is complete.\hfill{$\square$}

\medskip
\noindent
{\bf Remark 5.} It was proved by Pikhurko \cite{pikh01} that $v(P)\le c_d\cdot |P|$
holds for all $d$-dimensional lattice polytopes with nonempty interior lattice point, where $c_d$ is a constant depending only on $d$. Thus, by Theorem 1 one can deduce
$$\log \left(\sum_{j=1}^m v(d,j)\right) \gg m^{{d-1}\over {d+1}},$$
which was proved by B\'ar\'any \cite{bara08} by a different method.

\bigskip\noindent
{\bf Acknowledgement.} I am grateful to Prof. B\'ar\'any and the referees for their comments and suggestions,
which essentially improve the quality of this paper.

\bibliographystyle{amsplain}

\bigskip
\bigskip
\noindent C. Zong, School of Mathematical Sciences, Peking
University, Beijing 100871, China

\noindent E-mail: cmzong@math.pku.edu.cn

\end{document}